\documentclass[twoside,leqno,twocolumn,review]{article}
\usepackage{graphicx}

\usepackage{ltexpprt}
\usepackage{hyperref}
\usepackage{amssymb, amsmath} 
\usepackage{comment}
\usepackage{algorithm}
\usepackage{algpseudocode}
\usepackage{xcolor}

\DeclareMathOperator*{\argmin}{argmin}

\begin{document}
%\Crefname{ALC@unique}{Line}{Lines} % <- Preamble

\newcommand\relatedversion{}

\title{\Large Anderson Accelerated PMHSS for Complex-Symmetric Linear Systems \relatedversion}

\author{M\aa ns I. Andersson \thanks{KTH Royal Institute of Technology (mansande@kth.se)} \and Felix Liu \and Stefano Markidis} 

\date{}

\maketitle
% Copyright Statement
% When submitting your final paper to a SIAM proceedings, it is requested that you include
% the appropriate copyright in the footer of the paper.  The copyright added should be
% consistent with the copyright selected on the copyright form submitted with the paper.
% Please note that "20XX" should be changed to the year of the meeting.

% Default Copyright Statement 
%\fancyfoot[R]{\scriptsize{Copyright \textcopyright\ 2023 by SIAM\\
%Unauthorized reproduction of this article is prohibited}}

% Depending on which copyright you agree to when you sign the copyright form, the copyright
% can be changed to one of the following after commenting out the default copyright statement
% above.

\fancyfoot[R]{\scriptsize{Copyright \textcopyright\ 2023\\
Copyright for this paper is retained by authors}}

\pagenumbering{arabic}

\begin{abstract} \small\baselineskip=9pt 
This paper presents the design and development of an Anderson Accelerated Preconditioned Modified Hermitian and Skew-Hermitian Splitting (AA-PMHSS) method for solving complex-symmetric linear systems with application to electromagnetics problems, such as wave scattering and eddy currents. While it has been shown that the Anderson Acceleration of real linear systems is essentially equivalent to GMRES, we show here that the formulation using Anderson acceleration leads to a more performant method. We show relatively good robustness compared to existing preconditioned GMRES methods and significantly better performance due to the faster evaluation of the preconditioner. In particular, AA-PMHSS can be applied to solve problems and equations arising from electromagnetics, such as time-harmonic eddy current simulations discretized with the Finite Element Method. 
We also evaluate three test systems present in previous literature. We show that the method is competitive with two types of preconditioned GMRES. One of the significant advantages of these methods is that the convergence rate is independent of the discretization size. 

\end{abstract}

\section{Introduction}
\label{intro}
% preconditioned modified Hermitian and skew-Hermitian splitting = PMHSS

The solution of complex-valued linear systems arises in many scientific and engineering problems such as the time-dependent Schrödinger equation~\cite{li2018fast}, structural mechanics~\cite{adams2007algebraic}, Electromagnetism problems such as wave scattering \cite{baumann2018convergence} and time-harmonic eddy current problems \cite{bai2012block,rodriguez2010eddy}. Antenna problems solved with the method of moments (MoM) or other Boundary Element Methods often need to numerically calculate large complex-valued systems of linear equations~\cite{wang2007short}. 

A widely studied family of methods to solve these systems is the splitting technique, in particular the skew-Hermitian splitting (HSS). The basic idea of the HSS is based on the combination of operator splitting and alternating direction iteration (ADI) \cite{bai2007accelerated}. A special case of these type of solvers is called Preconditioned Modified Hermitian and Skew-Hermitian Splitting (PMHSS)~\cite{bai2011preconditioned}. PMHSS can be used as a solver for complex-valued linear systems by itself, but perhaps more practically, PMHSS iterations can also be used to precondition linear systems prior to solving with other methods, such as GMRES.

In this work, we propose to accelerate the PMHSS iteration method with the well-known Anderson Acceleration (AA)~\cite{anderson1965iterative} technique for improving the rate of convergence of fixed point iterations. This provides a different approach to solving complex symmetric linear systems based on PMHSS, which we will refer to as Anderson accelerated PMHSS (AA-PMHSS). In view of the known \emph{essential equivalence} of GMRES and Anderson acceleration in the linear case \cite{walker2011anderson, ni2009anderson}, the difference between AA-PMHSS and PMHSS-preconditioned GMRES may seem unclear. However, AA-PMHSS has the advantage of a by-design cheaper inner solve, when an iterative solver is used, which leads to faster solution times and fewer inner iterations in many cases. In this paper, we will discuss some theoretical as well as empirical results related to AA-PMHSS. In particular, we look at how it relates in terms of spectral properties with other similar methods, and we implement a prototype of AA-PMHSS in Matlab to investigate its performance and numerical properties empirically.

The main contributions of this work are the following:
\begin{enumerate}
    \item We propose a new method, based on AA and PMHSS, for accelerating solving symmetric complex-valued linear systems. We utilize the relationship between AA and GMRES to study the convergence of the method. 
    \item We design and implement a prototype version of the AA-PMHSS algorithm and assess its properties and performance in terms of iterations, computational efficiency, and wall-clock time.
    \item We compare the performance of the AA-PMHSS with relevant existing methods, including PMHSS preconditioned GMRES, PRESB preconditioned GMRES for solving relevant systems. We show that AA-PMHSS is a competitive method with current state-of-the-art preconditioners, and in most cases, superior.
\end{enumerate}

The remainder of the paper is organized as follows. We begin with the background in Section~\ref{background}, where we present the PMHSS iteration and some of its properties. Then we formulate Anderson acceleration of PMHSS and a short description of preconditioned GMRES with PMHSS and the C-to-R transform with PRESB preconditioner. Section~\ref{implementation} describes the implementation and in Section~\ref{exp_setup} we present the hardware and software setup and the test system. In Section~\ref{results}, we first show the convergence of the different methods on three different discretizations and for two different frequencies, followed by the convergence of the inner solver and finally a table describing wall-clock time, iterations, and the different errors. Section~\ref{rel_work} discusses the previous work on the solution of complex-valued linear systems, together with their applications. Finally, Section~\ref{conclusion} summarizes the main results and concludes the paper.

\section{Background}
\label{background}
We begin by presenting the background of our problem, and the necessary prerequisites to develop our Anderson accelerated version of PMHSS. Our focus is on problems that arise from discretized partial differential equations in electromagnetics, which can be both large-scale and have physical properties that lead to ill-conditioned systems. Consider the linear system
\begin{align}
\label{eq:Axb}(A+iB)\mathbf{x} = \mathbf{b},
\end{align}
where $A, B \in \mathbb{R}^{N \times N}$ and $\mathbf{x},  \mathbf{b} \in \mathbb{C}^{N}$.  In particular, we will focus on the cases where $A$ is symmetric positive definite (SPD) and $B$ symmetric positive semidefinite (SPSD), for which previous work's theory holds. 
\subsection{PMHSS Iteration}

PMHSS iteration, as described in \cite{bai2011preconditioned} and \cite{axelsson2014comparison}, has the form of an alternating direction type method. The iteration is given by the following relations:
\begin{align}
    (\alpha V + A) \hat{\mathbf{x}}_{k+1/2} &= (\alpha V - iB) \hat{\mathbf{x}}_k + \mathbf{b} \\
    (\alpha V + B) \hat{\mathbf{x}}_{k+1} &= (\alpha V + iA) \hat{\mathbf{x}}_{k+1/2} - i \mathbf{b},
\end{align}
and can be viewed as a type of fixed point iteration, which will be useful to us later. For $k = 0,1,2, ...$. For the particular choice of preconditioning parameter $\alpha = 1$ and preconditioner $V = A$ we have the following formulation used in \cite{axelsson2014comparison}:  
\begin{align}b
    (2A) \hat{\mathbf{x}}_{k+1/2} &= (A - iB) \hat{\mathbf{x}}_k + \mathbf{b} \\
    (A + B) \hat{\mathbf{x}}_{k+1} &= (1 + i)A \hat{\mathbf{x}}_{k+1/2} - i \mathbf{b}.
\end{align}
Substituting (2.4) into (2.5) gives us a single system per update, 
\begin{align}
  \label{eq:PMHSSit}  (A + B) \hat{\mathbf{x}}_{k+1} &= \frac{(1 + i)}{2} (A - iB) \hat{\mathbf{x}}_k + \frac{(1 - i)}{2} \mathbf{b},
\end{align}

It has been shown that the method has guaranteed convergence to a solution of \autoref{eq:Axb} independently of the initial guess, assuming that A is SPD and B is SPSD. Note that the left-hand side is transformed into a real matrix, but the right-hand side still requires a complex matrix-vector product. We will focus on the case when \autoref{eq:PMHSSit} is solved with an iterative scheme such as a Krylov space method, in our case the Conjugate Gradient (CG) method. %or a multigrid method. 
\subsection{Convergence of PMHSS iteration} 
The convergence of PMHSS iteration (under suitable assumptions) has, naturally, been analyzed previously. We show some simplified main points to highlight behavior used later. As discussed, given our choice of $\alpha$ and preconditioner, the iterations simplify to the form stated in \autoref{eq:PMHSSit}. Multiplying by $(A + B)^{-1}$ on both sides leads us naturally to define the following update matrix $\Psi$ and update vector $\mathbf{c}$:
\begin{align}
   \Psi &= \frac{(1 + i)}{2} (A + B)^{-1}(A - iB) \\
    \mathbf{c} &= \frac{(1 - i)}{2}(A + B)^{-1}\mathbf{b}. 
\end{align}
Assuming an initial guess $\mathbf{\hat{x}_0} = \mathbf{0}$ (and thus $\mathbf{\hat{x}_1} = \mathbf{c}$), repeated application of the update rule in \autoref{eq:PMHSSit} gives us the following relation: 
\begin{align}
\label{eq:update}
    \mathbf{\hat{x}}_{k+1} =& \Psi(...(\Psi \mathbf{c} + \mathbf{c})+...)+\mathbf{c} \\
    =& (\Psi^k + ... + \Psi + I)\mathbf{c}.
\end{align}
The expression in the parentheses is a geometric sum (of matrices), which has the closed form $(I - \Psi)^{-1}(I - \Psi^k)$. Hence, we have 
\begin{align}
\label{eq:inner}
\mathbf{\hat{x}_{k}} = (I - \Psi)^{-1}(I - \Psi^{k})\mathbf{c}.
\end{align}
When the spectral radius $\rho(\Psi) < 1$ as in \cite{bai2010modified,bai2011preconditioned}, we have that $(I - \Psi)^{-1}(I - \Psi^k) \rightarrow (I - \Psi)^{-1}$ as $k \rightarrow \infty$. Thus, in the limit, the solution obtained satisfies:
\begin{align}
\label{eq:finaleq}
(I - \Psi)\mathbf{\hat{x}} = \mathbf{c},
\end{align}
which a straightforward computation shows is equivalent to \autoref{eq:Axb}.

\subsection{Inner solver}
\label{sec:inner}

The previous system \autoref{eq:PMHSSit} is solved for each step of the PMHSS iteration. As the iteration progresses, we get a better initial guess for the solution of the next system by re-using $\mathbf{x}_k$. 

For example, if we use CG iteration to solve the inner system. 
Assuming that we want to calculate the $\mathbf{x}_{k+1}^*$, the true solution at $k+1$, we have the following estimated initial error $e_0$ for the inner linear solver:
\begin{align}
   \ ||e_{0}(k)|| &= ||\mathbf{x}^*_{k+1}-\mathbf{\hat{x}}_{k}||\\
            &\approx ||\mathbf{\hat{x}}_{k+1}-\mathbf{\hat{x}}_{k}||\\
           \label{eq:e0inner} 
           &= ||(I - \Psi)^{-1}(\Psi^{k} - \Psi^{k+1})\mathbf{c}||\\
           &= ||(I - \Psi)^{-1}(I - \Psi)\Psi^{k}\mathbf{c}||\\
            &\leq ||\Psi^{k} || \cdot ||\mathbf{c}||. \label{eq:initerror}
\end{align}
Where $||\Psi^{k} ||$ will decrease with $k$. We realize here that the inner and outer iteration is naturally converging to the same solution. 
\begin{comment}
In comparison, using $\mathbf{x}_0 = \mathbf{0}?$
\begin{align}
\label{eq:e0inner0} 
    ||e_{0}(\mathbf{0})|| &= ||\mathbf{x}^*_{k+1}-\mathbf{0}||\\
            &\approx ||\mathbf{\hat{x}}_{k+1}||\\
            &= ||(I - \Psi)^{-1}(I - \Psi^{k+1})\mathbf{c}|| \\
            &\leq ||(I - \Psi)^{-1}|| \cdot ||(I - \Psi^{k+1}) || \cdot ||\mathbf{c}||
\end{align}
Which is \textit{k}-independent and clearly larger than \autoref{eq:initerror}:
\begin{align}
    ||e_{0}|| \approx || \Psi ||^k \cdot ||e_{0}(\mathbf{0})||.
\end{align}
\end{comment}

\subsection{Anderson Acceleration}
Anderson acceleration~\cite{anderson1965iterative} (known as Anderson mixing in some fields) is a method for accelerating fixed-point iteration schemes for solving problems of the form $\mathbf{x} = f(\mathbf{x})$. A simple method for finding fixed points is to simply compute $\mathbf{x}_{i+1} = f(\mathbf{x}_i)$ until convergence by some criterion, however convergence of such a method can be slow. The idea behind Anderson acceleration is to search for the fixed point of $f(\mathbf{x}_k) = \mathbf{x}_{k+1} := \mathbf{f}_k$ by restating it as an optimization problem: $g(\mathbf{x}_k) = f(\mathbf{x}_k) - \mathbf{x}_{k} := \mathbf{g}_k$ and find $\mathbf{x}$ such that $g(\mathbf{x}) = \mathbf{0}$. From a least-squares perspective, this leads naturally to consider the problem of minimizing $||g(\mathbf{x})||_2$. Anderson acceleration considers least squares solutions in the space spanned by the previous $\mathbf{g}_k$:
\begin{equation}
    \argmin_{\alpha_k \in A_k} ||G_k \alpha||_2, \\ \quad A_k = \{\mathbf{\alpha} \in \mathbb{R}^k : \sum_i \alpha_i = 1 \},
\end{equation}
where $G_k = $[$\mathbf{g}_0$, ..., $\mathbf{g}_k$]. Then the current best solution is updated accordingly. 
\begin{equation}
    \mathbf{x} = F_k \mathbf{\alpha}_k \quad F_k = [f_0, ... , f_k]
\end{equation}
There are many similarities with GMRES that are described in the next section, with the main difference being that we find the solution as a minimization of $r(\mathbf{x})$ instead of $g(\mathbf{x})$.

The suggested algorithm AA-PMHSS is described below and is based on a reformulation of AA found in \cite{walker2011anderson,fang2009two}, which simplifies the implementation of the optimization problem by making it unconstrained. Formulating a truncated version where only the latest $m$ iterations are used, this is sometimes called a \textit{windowed} AA and can reduce the memory needed for the method, at the cost of a reduced convergence rate. We will not use truncation in the remainder. It is also possible to write a relaxed formulation; however, we omit that in this study. 

\begin{algorithm}
\caption{AA-PMHSS iteration, assuming using a CG solver for the inner system in 1 and a QR direct linear solver for 3}
\label{alg:aa-pmhss}
%\SetAlgoLined
%\KwResult{Write here the result }

\begin{algorithmic}[1]
\State Solve : $(A+B)\mathbf{x}_{k+1}=\frac{(1+i)}{2}(A-iB)\mathbf{\bar{x}}_k + \frac{(1-i)}{2}\mathbf{b}$
\State $X_k = $[$\mathbf{x}_1-\mathbf{x}_0, ... , \mathbf{x}_{k-1}-\mathbf{x}_k$]$ \newline G_k = $[$\mathbf{g}_1-\mathbf{g}_0, ... , \mathbf{g}_{k-1}-\mathbf{g}_k$]$  \newline \mathbf{g}_k$ = $\mathbf{x}_{k+1}-  \mathbf{\bar{x}}_k$
 \State Solve : $\argmin\limits_{\alpha_k \in A_k} ||\mathbf{g}_k - G_k\mathbf{\alpha}||_2$, \quad $A_k$ = $\{\alpha \in \mathbb{C}^k \}$ %: \sum_i \alpha_i = 1 \}$\\
\State$\mathbf{\bar{x}}_{k+1} = \mathbf{x}_k  +  \mathbf{g}_k - (X_k+G_k)\mathbf{\alpha}_k$
\end{algorithmic}
\end{algorithm}
\subsection{Preconditioned GMRES}
The generalized minimal residual method (GMRES) \cite{saad1986gmres} is a Krylov subspace method for the solution of linear systems of equations $C\mathbf{x} = \mathbf{b}$, which can handle both non-symmetric and indefinite matrices. The key idea, as in all Krylov methods, lies in searching for solutions in a Krylov subspace of the linear system, which has the following form:
\begin{align}
    \mathcal{K}(C,\mathbf{v}) := \text{span}(C\mathbf{v},C^2\mathbf{v},C^3\mathbf{v},..). 
\end{align}

 The idea in GMRES is to look for least squares solutions (in the Euclidean norm) to $C\mathbf{x} = \mathbf{b}$ for $\mathbf{x} \in \mathcal{K}(C,\mathbf{v})$. The implementation in practice relies on Arnoldi's method to iteratively build an orthonormal basis of $\mathcal{K}(C,\mathbf{v})$. Due to the way GMRES is formulated, the solution of the resulting linear least squares problem can be carried out efficiently using QR-factorization. The QR-factorization itself can also be updated cheaply throughout the iterations and does not to be re-formed.

For $C=A+iB$ and $\mathbf{v}=\mathbf{r}_0$, GMRES finds the minimal residual to Equation \ref{eq:Axb} in a Krylov subspace, that is the $m$\textsuperscript{th} iteration satisfies
\begin{align}
    \mathbf{x}_m = \argmin_{x\in K_m(C,\mathbf{r}_0)} || r(\mathbf{x}) ||_2.
\end{align}

The proposed preconditioner from previous work is presented in Algorithm \ref{alg:precPMHSS}. The PMHSS-GMRES was implemented as in~\cite{axelsson2014comparison} with Matlab's built-in GMRES (which is based on \cite{walker1988implementation}). 

\begin{algorithm}
\caption{PMHSS precondition, where $\mathbf{q}$ is the current residual in the iterative method. Note that the second line is not strictly necessary.}
\label{alg:precPMHSS}
\begin{algorithmic}[1]
\State Solve $(A+B)\mathbf{z}=\mathbf{q}(\mathbf{x})$
\State Set $\mathbf{x}= \frac{(1-i)}{2}\mathbf{z}$
\end{algorithmic}
\end{algorithm}
We can perform a brief analysis of the eigenvalue distribution of the preconditioned system. Using Algorithm \ref{alg:precPMHSS} assuming that A is SPD and B is SPSD so that the eigenvalues $\mu = \mu(A^{-1}B)$ are real and non-negative, we have that $(A+B)^{-1}(A+iB) = (I + (i-1)(I +A^{-1}B)^{-1}(A^{-1}B))$ thus the preconditioned eigenvalues are:
\begin{align}
\label{eq:spectrum-pmhss-gmres}
    \lambda_{\text{PMHSS-GMRES}} = 1 + (i - 1)\frac{\mu}{\mu+1}.
\end{align}
We then have, if $\mu_{\text{min}} \leq \mu_{} \leq \mu_{\text{max}}$,
\begin{align}
 \frac{1}{1+\mu_{\text{max}}} \leq \mathfrak{Re}(\lambda) \leq \frac{1}{1+\mu_{\text{min}}}, \\
 \frac{\mu_{\text{min}}}{1+\mu_{\text{min}}} \leq \mathfrak{Im}(\lambda) \leq \frac{\mu_{\text{max}}}{1+\mu_{\text{max}}}.
\end{align}
We clearly see that $\lambda$ is not dependent on discretization. %We can note that $\lambda$ still generally is complex. 

Note that the preconditioned system in \text{Algorithm \ref{alg:precPMHSS}} is not equivalent with \autoref{eq:PMHSSit} as it differs with a factor $-i$. Therefore, we cannot reuse the current best solution as a guess for the inner solver in PMHSS-GMRES, which can be done in the PMHSS iteration and AA-PMHSS.

\subsection{Relationship between Anderson Acceleration and GMRES}

Multiple works have shown that the AA solution of a linear system and GMRES is \textit{essentially equivalent} \cite{ni2009anderson,walker2011anderson,de2021anderson}. More precisely this means that, at any k+1 iteration of AA with a solution $x_{\text{AA}}$, one can construct the solution $x_{\text{GMRES}}$. We must mention that no proof was found to hold for the complex-valued case. The non-trivial assumption used is that $\sum_i \alpha_i = 1$. This does not hold generally when $\alpha \in \mathbb{C}^{k}$ which hinders us from reusing previous proofs. Disregarding that, we analyze the convergence of the AA-PMHSS method by restating the GMRES equivalent. Let
\begin{align}
    f(\mathbf{x}) = C\mathbf{x} + \mathbf{c}.
\end{align}
Solving for a fixed-point with AA is \textit{essentially} equivalent to solving the following system with GMRES:
\begin{align}
(I - C)\mathbf{\hat{x}} = \mathbf{c},
\end{align}
which we recognize from \autoref{eq:finaleq}. 
Since every GMRES and AA iteration have corresponding solutions, we can model the convergence of AA by using GMRES. 

In our case for AA-PMHSS we have the update equation \autoref{eq:PMHSSit} which, when restated, gives the following GMRES systems, which is the same as \autoref{eq:finaleq}:
\begin{align}
(1+i)(I-\frac{(1+i)}{2}(A+B)^{-1}(A-iB))\mathbf{x} = \mathbf{b}.
\end{align}
Which, by similar analysis as before, gives us the preconditioned eigenvalues for the AA-PMHSS iteration as
\begin{align}
    \lambda_{\text{AA-PMHSS}} &= 1+i-i\Bar{\lambda}_{\text{PMHSS-GMRES}} \\
    &= 1 + (i - 1)\frac{\mu}{\mu+1}. %\lambda_{PMHSS-GMRES(C)}
\end{align}
This is the same spectrum as PMHSS-GMRES found in \autoref{eq:spectrum-pmhss-gmres}. As such, we can assume that the AA-PMHSS, if stable, will converge in as many iterations as the PMHSS-GMRES. To summarize the section:
\begin{itemize}
    \item The convergence of AA-PMHSS will be similar to PMHSS-GMRES and discretization independent, as the eigenvalues will be contained in a domain of the same size. Since the GMRES equivalent of AA-PMHSS has the same expected convergence as PMHSS-GMRES. 
    \item The inner system (the PMHSS iteration \autoref{eq:PMHSSit}) in AA-PMHSS is approaching the same system as the outer iterations \autoref{eq:Axb}. Therefore, the number of inner iterations will decrease linearly as the outer iteration progresses. 
    \item We can expect that the reduction in cost for the inner iteration will decrease faster in AA-PMHSS than in the standard PMHSS iteration due to the faster convergence of the outer iteration. 
\end{itemize}
%\autoref{eq:inner}

\subsection{C-to-R transform}
A common method to solve complex-valued systems is to use the complex-to-real transform (C-to-R). The system \autoref{eq:Axb} can be written as a real system on the form, 
\begin{align}
    \begin{bmatrix}
    A & -B \\
    B & A 
    \end{bmatrix} 
    \begin{bmatrix}
    \mathbf{x} \\
    \mathbf{y}
    \end{bmatrix}
    = 
    \begin{bmatrix}
    \mathbf{a} \\
    \mathbf{b}
    \end{bmatrix}
\end{align}
There are multiple rewrites that can be situationally more useful. This system can be solved efficiently with the AGMG multigrid solver or with a PRESB preconditioned GMRES described in the next section. Which, we will focus on as it shares PMHSS' discretization independent convergence. 
\subsection{PRESB}
The C-to-R representation enables another interesting method, the PREconditioned Square Block (PRESB) presented in \cite{axelsson2014comparison,AXELSSON2020286,axelsson2000real}. The PRESB preconditioner is given in (2.30) and admits the following factorization.
\begin{align}
\mathcal{P}_{\text{PRESB}} &= \begin{bmatrix}
A & B \\
B & A+2B
\end{bmatrix} \\
&= \begin{bmatrix}
I & -I \\
0 & I
\end{bmatrix}
\begin{bmatrix}
A+B & 0 \\
B & A+B
\end{bmatrix}
 \begin{bmatrix}
I & I \\
0 & I
\end{bmatrix}
\end{align}
Thus, the PRESB matrix is \emph{similar} to the middle matrix in (2.31). The spectral properties of our preconditioner can thus be analyzed by considering the system:
\begin{align}
\mathcal{P}_{\text{PRESB}}^{-1}\mathcal{A} = \begin{bmatrix}
A+B & 0 \\
B & A+B
\end{bmatrix}^{-1}\begin{bmatrix}
A & -B \\
B & A
\end{bmatrix}.
\end{align}

From \cite{axelsson2014comparison}, we find the following spectral bounds, assuming the same conditions as with the PMHSS analysis. 
\begin{align}
    \lambda_{\text{min}} = \frac{1}{2}, \quad \lambda_{\text{max}} = 1 
\end{align}
The implementation of the algorithm is described in \text{Algorithm \ref{alg:precPRESB}}, and notably consists of two linear systems to solve internally.
\begin{algorithm}
\caption{PRESB preconditioner \cite{axelsson2014comparison}, where $\mathbf{p}$ and $\mathbf{q}$ are the current real and imaginary residuals.}
%\SetAlgoLined
%\KwResult{Write here the result }
\label{alg:precPRESB}
\begin{algorithmic}[1]
\State Let $H_1 = A+B_1$, $H_2 = A+B_2$
\State Solve $H_1 \mathbf{h}= \mathbf{p} + \mathbf{q}$
\State Solve $H_2 \mathbf{y} = \mathbf{q} - B_1 \mathbf{h}$
\State Compute $\mathbf{x} = \mathbf{h} - \mathbf{y}$
\end{algorithmic}
\end{algorithm}
\section{Implementation}
\label{implementation}
The implementations of C-to-R and PHMSS preconditioned GMRES are done using Matlab's GMRES function, the full AA-PMHSS algorithm written in Matlab using the built-in multithreaded routines for linear algebra. The inner solver is by default set to Conjugate Gradient (CG), with the possibility of further preconditioning. For some problems, other inner solvers may be needed. The minimization problem in AA-PMHSS is solved with Matlab's QR factorization and the backslash operator.  We also use the Matlab interfaced version of AGMG 3.3.5~\cite{notay2010aggregation,notay2012aggregation,napov2012algebraic}.

\section{Experimental Setup}
\label{exp_setup}
We now describe the hardware and software environment and then the test case. All timing experiments were performed using MATLAB R2022b on a local workstation with an AMD Ryzen 9 7900x CPU with 64 GB of DDR5 DRAM with all types of energy saving and clock boost disabled. The timings are performed with Matlab's \texttt{tic} and \texttt{toc} commands, ten runs are performed, and the standard deviation is collected and referred to as $\sigma$.

\subsubsection*{Eddy Current}
The following problem is a time-harmonic eddy current with a conditioning parameter $\omega$. The system matrix comes from a high-order Nédléc finite element method. The studied frequency governs the conditioning of the system matrix. We have two cases, high (1000) and low (0.0001). The system is complex-symmetric and positive definite. We run three discretization sizes N = 293, 2903, and 25602 to investigate the size-independent convergence rate. We generate a random complex valued right-hand side, which is reused for all methods, and use the zero vector as an initial guess. If not stated we have a relative tolerance set to $10^{-8}$ for the full solver calculated with $r_c(x) = ||(A+iB)x-c||/||c||$ and $10^{-12}$ for the inner solver. For the AA-PMHSS method, we use $||g(x^*)||$ as the stopping criteria for timings, but we evaluate the error explicitly for the complete results. 

The three following problems are finite differences discretizations found in the literature defined on the unit square and are used to generate large benchmark problems. 
\subsubsection*{Padé approximation} 
Previous references \cite{axelsson2000real,bai2010modified} have the \autoref{eq:Axb} on the form, 
\begin{align}
\left[\left(L + \frac{3-\sqrt{3}}{h}I\right) + i \left(L + \frac{3+\sqrt{3}}{h}I\right)\right] \mathbf{x} = \mathbf{b}    
\end{align}
Where $h$ is the spatial discretization length, $I$ is the identity, and $L$ is a 5-point finite-difference discretization of the negative Laplace operator on the unit square. This complex-valued linear system appears in the approximated time integration of parabolic PDEs \cite{axelsson2000real}. We use the right-hand side from the same source. 
\subsubsection*{Inhomogeneous Helmholtz equation or Shifted $\omega$-system} The Helmholtz equation appears in \textit{acoustics} and \textit{electromagnetism}. 
\begin{align}
    \bigg[L + (\mu + i \omega) I\bigg] \mathbf{x}  = \mathbf{b}.
\end{align}
Here we use the coefficients $\mu = 0$ and $\omega = 0.01$ found \cite{axelsson2000real,axelsson2014comparison}. This particular system appears when solving Ordinary Differential Equations on the form,
\begin{align}
\dot{\mathbf{q}} + L\mathbf{q(t)} = f(t),
\end{align}
by assuming a periodic \textit{ansatz} for $\mathbf{q}$. We use a right-hand side with uniformly distributed values between negative one and one for the real and imaginary entries.
\subsubsection*{Equation of motion} 
Consider the linear system on the form
\begin{align}
    \bigg[ \left( K - \omega^2 M \right) + i \left(\omega C_V + C_H \right) \bigg] \mathbf{x} = \mathbf{b}.
\end{align}
Where M and K are the inertia- and the stiffness matrices, we have $C_V$ and $C_H$ are the viscus- and hysteretic dampening matrices. From previous work \cite{axelsson2000real,bai2010modified,wu2017modified} we use $K=L$, $M=I$, $C_v=10I$, and $C_H=\mu K$. The scalars $\omega$ and $\mu$ are set to $\pi$ and $0.02$ respectively. We use the same right-hand side as with the shifted-$\omega$ benchmark. This system arises when solving the direct frequency response to mechanical systems, which can be written on the form, 
%systems on the form in the frequency domain,
%the direct frequency response mechanical systems 
\begin{align}
    M\ddot{\mathbf{q}} + (C_V + \frac{1}{\omega} C_H)\dot{\mathbf{q}} + K \mathbf{q} = \mathbf{p}.
\end{align}
For more information \cite{benzi2008block,feriani2000iterative}.
\section{Results}
\label{results}
In this section, we provide some numerical results to compare the performance of our proposed method (AA-PMHSS) with PMHSS and PRESB preconditioned GMRES, as well as unpreconditioned GMRES. Our goal is to investigate the practical implications of the theoretical discussion in \autoref{background}, as well as to compare the performance of AA-PMHSS with other related methods. We begin by showing the convergence rate of the solvers on three different sizes of systems to study how the rate of convergence of the examined methods depends on problem size. Then, for the solvers preconditioned with PMHSS (AA-PMHSS and PMHSS-GMRES), we look at how the solver's convergence depends on the accuracy and tolerances used for the inner solver. Lastly, we collect more detailed performance results, including the number of outer and inner iterations, tolerances, and wall time on the test problems described in \autoref{exp_setup}.

\begin{figure}[ht]
    \centering
    \includegraphics[width=0.5\textwidth]{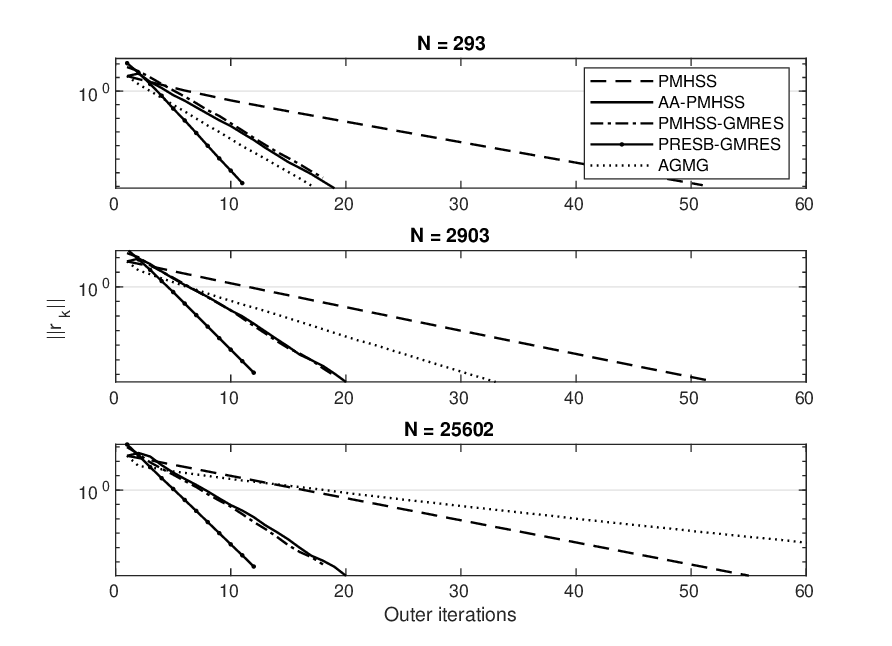}
    \caption{Convergence of three discretization sizes, showing that the upper bound of the number of outer iterations needed is independent of the number of degrees of freedom. AGMG is included as a reference as a standard linear solver.}
    \label{fig:discretizationandconvergence}
\end{figure}

\autoref{fig:discretizationandconvergence} shows the convergence in the residual norm of the different solvers as a function of the outer iterations for three different linear system sizes. In our results, we see that the convergence of the three PMHSS-based methods (two with PMHSS-preconditioning and pure PMHSS iteration) as well as PRESB-GMRES is essentially independent of the system size, a property which does not hold for the general purpose solver, algebraic-multigrid based AGMG, which is included for reference. This size-independence is what one would expect from the spectral bounds for the PMHSS-based methods in \autoref{background}. 
 
\begin{figure*}[ht]      
    \centering
    \includegraphics[width=0.50\textwidth]{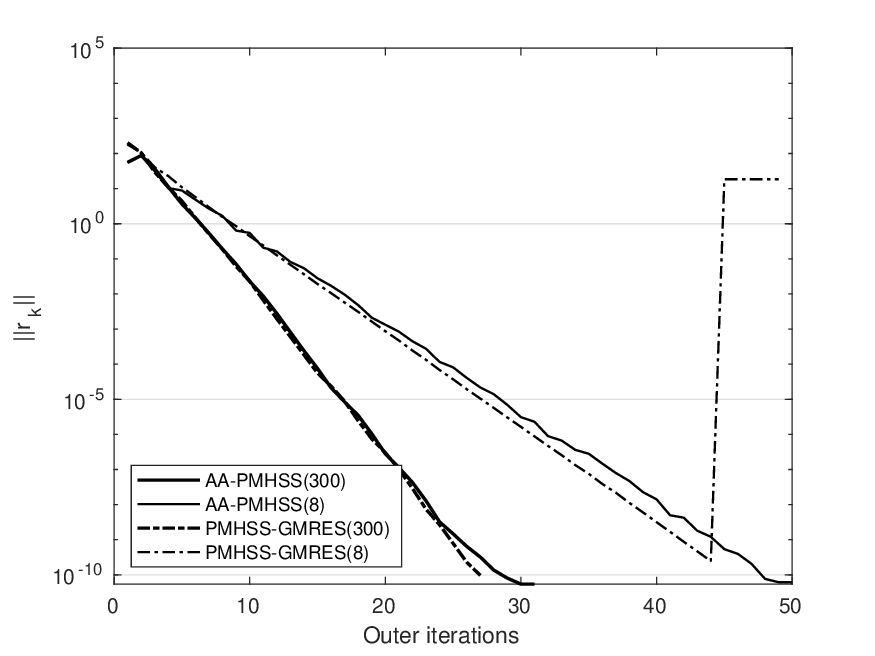}\includegraphics[width=0.50\textwidth]{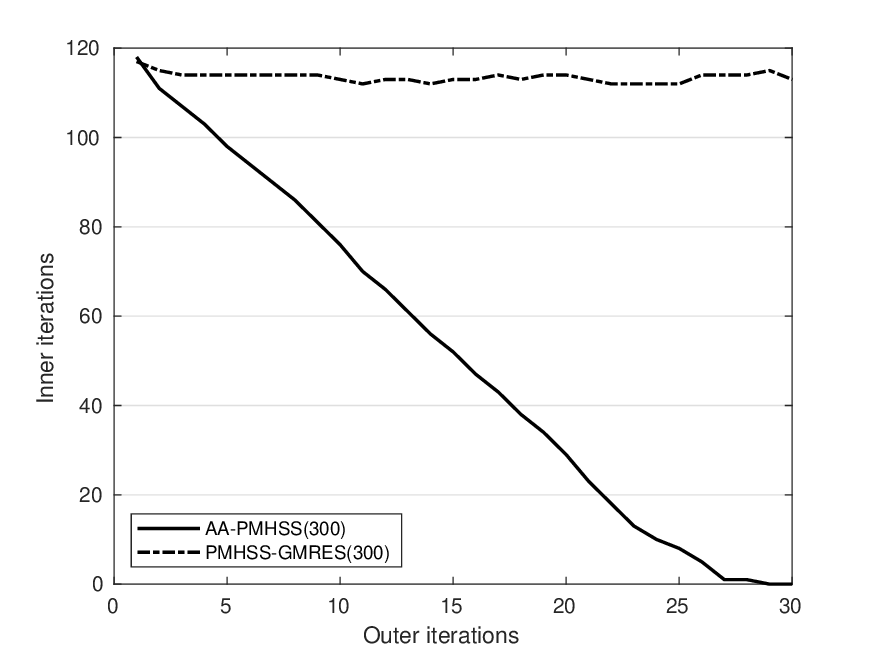}
    \caption{\textbf{(left)} The convergence of AA-PMHSS and PMHSS-GMRES with a fixed maximal number of iterations for the inner CG-solver (300) and (8). This run was performed on the medium-sized eddy current problem and the tolerance was set to 1e-12. \textbf{(right)} The number of inner iterations needed for the inner solver to converge to machine precision for the solvers AA-PMHSS and PMHSS-GMRES. This run was performed on the medium-sized problem for $\omega = 0.0001$.  For AA-PMHSS, we find that the current best solution provides a good initial guess for the inner solver. PMHSS-GMRES used the zero vector as initial guess.}\label{fig:maxinner} 
\label{fig:nrinnerit}
\end{figure*}

Furthermore, the rate of convergence observed in our numerical results follows what one would expect from the theoretical spectral bounds discussed in \autoref{background}. PRESB-GMRES shows the fastest convergence in terms of outer iterations in \autoref{fig:discretizationandconvergence} and also has the tightest spectral bounds. Furthermore, AA-PMHSS and PMHSS-GMRES show almost the same convergence rate in the numerical experiments, which is also what one would expect from the theoretical analysis. This may again be unsurprising given the aforementioned essential equivalence between Anderson Acceleration and GMRES discussed in \cite{ni2009anderson,walker2011anderson}.
The left figure in \autoref{fig:maxinner} shows the convergence, again in the residual norm, as a function of outer iterations for AA-PMHSS and PMHSS-GMRES, but where we set the maximum number of inner iterations to 300 and 8. As expected, the convergence in terms of outer iterations is worsened when we use a smaller number of inner iterations, but the convergence to low residual norms is retained for all cases but PMHSS-GMRES with 8 inner iterations. This trade-off between accuracy of the inner solver and convergence of the outer iterations might be more than worthwhile in terms of total time to solution. The convergence of the AA-PMHSS solver starts to deteriorate at low residual norms, which PMHSS-GMRES (in the 300 inner iteration case) does not. One explanation for this observation can be seen in the right figure (in \autoref{fig:maxinner}), where the number of inner iterations in AA-PMHSS becomes zero in later outer iterations, indicating that the outer solver's convergence breaks down at that point.
The right plot in \autoref{fig:nrinnerit} shows the number of inner CG iterations required \emph{in each} outer iteration to converge to machine precision for AA-PMHSS and PMHSS-GMRES. The first observation one can make is that 300 was a pessimistic upper bound on the inner iterations in this case, as both methods manage to solve to the required tolerance in significantly fewer iterations. The more encouraging result is that the number of inner iterations required for our AA-PMHSS method decreases (approximately linearly) as the outer iterations progress. The reason for this is our ability to use the current best solution as the starting point for the inner solver in the next iteration; this is discussed in greater detail in \autoref{sec:inner}. 

\begin{comment}
\begin{figure}[h!]
    \centering
    \includegraphics[width=0.5\textwidth]{Figures/nrinnerit.eps}
    \caption{The number of inner iterations needed for the inner solver to converge to machine precision for the solvers AA-PMHSS and PMHSS-GMRES. This run was performed on the medium-sized problem for $\omega = 0.0001$.  For AA-PMHSS, we find that the current best solution provides a good initial guess for the inner solver. PMHSS-GMRES used the zero vector as initial guess.}
   \label{fig:nrinnerit}
\end{figure}
\end{comment}

We conclude the results by looking at the more comprehensive results shown in Tables~\ref{table:AA-PHMS} and \ref{table:perf}. Table~\ref{table:AA-PHMS} shows a performance comparison on the largest (N=25602) eddy current system. For AA-PMHSS, we also show results where the maximum number of CG iterations have been capped at 50 and 25 respectively. Note also that the results in \autoref{table:AA-PHMS} include tests where a direct linear solver (LU) is used for the inner solve, for comparison.

Looking at the number of outer iterations for convergence for the different methods in \autoref{table:AA-PHMS} shows some interesting results. Firstly, we note that the use of an iterative inner solver has not affected the convergence of the method in terms of outer iterations, almost at all, with the CG inner solver showing significantly faster solution times as well. The only exception is AA-PMHSS for the $\omega = 10^{-4}$ case, where the difference is only one outer iteration, as well as all the cases where the number of CG iterations in the inner solver is capped. While the capping of the inner CG iterations does affect the convergence in outer iterations (especially for the $\omega = 10^{-4}$ case), we see that it is worthwhile in terms of the \emph{total} number of inner iterations and solution time in this case. Finally, the convergence of the different methods in terms of outer iterations follows what one would expect from the theoretical discussion. AA-PMHSS and PMHSS-GMRES converge in a similar number of outer iterations for both problems, which aligns with the discussion regarding the essential equivalence of the two methods. PRESB-GMRES converges in the fewest number of outer iterations, which is expected given the tight spectral bounds, and pure PMHSS iteration has the slowest convergence.  

\autoref{table:perf} shows performance results for the Padé, shifted-$\omega$ and equation of motion matrices, for sizes $N=10000,40000,90000$ (square matrices). The best performing methods, in terms of runtime, are highlighted in bold. Overall, we see that AA-PMHSS performs well on most systems, being the \emph{fastest} for the Padé and equation of motion. Again, this can be attributed to requiring fewer inner iterations than PMHSS-GMRES, which is second in performance for those cases. PRESB-GMRES is the fastest on the shifted-$\omega$ systems, which is likely due to converging in much fewer outer iterations than the other methods; the PMHSS methods did not achieve the maximum number of outer iterations (that is, the size-independent convergence is not seen yet). Overall, PRESB-GMRES seems to have an advantage in converging in fewer outer iterations. However, it comes at the cost of each outer iteration being more expensive than the other methods and requiring two systems to be solved. In \autoref{table:perf_m50}, we again cap the inner iterations for the preconditioned methods but for the shifted-$\omega$ case. We find that the AA-PMHSS has a reduced convergence rate but performs almost a third as many inner iterations in total and cuts the wall-clock time in half. In contrast, the GMRES methods stagnate long before the tolerance. The PMHSS iteration is the second fastest with the capped inner iterations.  

Unsurprisingly, pure PMHSS iteration and unpreconditioned GMRES are the worst performing, and we can safely conclude that all the discussed preconditioning and acceleration strategies provide significant benefits in performance.

\begin{table*}[h]
\centering
\caption{Comparison of the different methods for solving the largest of the eddy current problems for both low-$\omega$ and high-$\omega$. With a direct linear solver (LU) and conjugate gradient (CG(n)) for solving the inner system. Where n refers to the maximum number of iterations allowed for the inner solver. The inner solvers' tolerance is set to $10^{-12}$. The reference solution $x^*$ was calculated with Matlab's \texttt{backslash}.\\
} 
\label{table:AA-PHMS}
\begin{scriptsize}
\begin{tabular}{lccccccc}
Method & outer iter. & inner iter. & T ± $\sigma$ [s] & $||r||$ & $||x^*-x||$ & $||(A+B)^{-1}r||$ \\ 
\hline \\
$\omega = 10^{4}$ \\
\hline PMHSS & & & & & & & \\ 
\quad LU & 67 & - & 25.43 ± 0.592\% & 2.73197e-09 & 1.03058e-09 & 1.02994e-09 \\
\quad CG & 67 & 2491 & 1.162 ± 2.973\% & 2.73819e-09 & 1.03532e-09 & 1.03377e-09 \\
\hline AA-PMHSS & & & & & & & \\
\quad LU & 22 & - & 8.733 ± 0.4441\% & 1.71475e-09 & 6.70805e-10 & 5.38084e-10 \\
\quad CG & 22 & 824 & 0.5365 ± 0.447\% & 1.93437e-09 & 7.51107e-10 & 6.1096e-10 \\
\quad CG(50) & 23 & 765 & 0.514 ± 0.3678\% & 9.89885e-10 & 3.953e-10 & 3.25222e-10 \\
\quad CG(25) & 23 & 496 & 0.3881 ± 0.6916\% & 1.66947e-09 & 6.44053e-10 & 5.23233e-10 \\
\hline PMHSS-GMRES & & & & & & & \\ 
\quad LU & 23 & - & 9.009 ± 0.03202\% & 1.00433e-09 & 3.90251e-10 & 2.1615e-10 \\
\quad CG & 23 & 1679 & 0.8161 ± 0.0251\% & 1.01339e-09 & 3.9347e-10 & 3.05682e-10 \\
\hline PRESB-GMRES & & & & & & & \\
\quad LU & 14 & - & 10.95 ± 0.06775\% & 1.88122e-09 & 9.2189e-10 & 5.98921e-10 \\
\quad CG & 14 & 2262 & 0.6489 ± 0.01114\% & 1.8851e-09 & 9.22169e-10 & 5.99389e-10 \\
\hline \\ 
\\
$\omega = 10^{-4}$ \\
\hline PMHSS & & & & & & &  \\
\quad LU & 80 & - & 30.25 ± 0.6644\% & 5.62865e-11 & 1.13009e-09 & 1.12844e-09 \\
\quad CG & 78 & 12506 & 5.491 ± 1.285\% & 2.06198e-10 & 3.7422e-09 & 3.79037e-09 \\
\hline AA-PMHSS & & & & & & &  \\
\quad LU & 25 & - & 9.861 ± 0.1814\% & 3.58754e-10 & 6.31284e-10 & 5.02025e-10 \\
\quad CG & 26 & 4061 & 2.099 ± 0.5552\% & 5.69164e-10 & 4.03492e-09 & 3.88251e-09 \\
\quad CG(50) & 32 & 1489 & 0.9801 ± 1.043\% & 5.93055e-10 & 4.06493e-09 & 3.82198e-09 \\
\quad CG(25) & 43 & 1012 & 0.9603 ± 1.684\% & 2.29485e-10 & 5.43464e-09 & 5.05516e-09 \\
\hline PMHSS-GMRES & & & & & & &  \\
\quad LU & 23 & - & 9.572 ± 0.08122\% & 4.73915e-09 & 1.18588e-08 & 6.48492e-09 \\
\quad CG & 23 & 7500 & 3.33 ± 0.03415\% & 5.09979e-09 & 1.35183e-08 & 1.11445e-08 \\
\hline PRESB-GMRES & & & & & & &  \\
\quad LU & 15 & - & 11.67 ± 0.07115\% & 3.23933e-09 & 6.54481e-09 & 4.27435e-09 \\
\quad CG & 15 & 10184 & 2.686 ± 0.02967\% & 4.00906e-09 & 8.49476e-09 & 7.52229e-09 \\
\hline
\end{tabular}
\end{scriptsize}
\end{table*}

\begin{table*}[]
\centering
\caption{The stencil-based benchmarks with CG for solving the inner system. The inner solvers' tolerance is set to $10^{-12}$, and the maximum number of iterations is the problem size. The reference solution $x^*$ is calculated with Matlab's \texttt{backslash}.\\
} 
\label{table:perf}
\begin{scriptsize}
\begin{tabular}{lccccccc}
Method/size & outer iter. & inner iter. & T ± $\sigma$ [s] & $||r||$ & $||x^*-x||$ & $||(A+B)^{-1}r||$ \\ 
\hline \\
& & & \underbar{\textbf{Padé}} \\
PMHSS \\
\hline 
\quad 10000 & 33 & 4172 & 0.677228 ± 5.79\% & 3.44893e-16 & 6.24067e-15 & 4.99637e-15 \\
\quad 40000 & 34 & 6134 & 2.63969 ± 0.8817\% & 1.71351e-17 & 6.20946e-16 & 4.99868e-16 \\
\quad 90000 & 34 & 7531 & 6.13483 ± 0.3547\% & 4.10001e-18 & 2.16516e-16 & 1.74593e-16 \\
GMRES & & & & & & &\\
\hline
\quad 10000 & 133 & - & 0.638218 ± 1.767\% & 7.715e-15 & 6.87689e-14 & - \\
\quad 40000 & 188 & - & 3.4308 ± 1.579\% & 9.54629e-16 & 2.35029e-14 & - \\
\quad 90000 & 227 & - & 9.44645 ± 2.2\% & 2.9451e-16 & 1.0854e-14 & - \\
AA-PMHSS & & & & & & & \\
\hline 
\quad 10000 & 10 & 1300 & \textbf{0.238979} ± 5.973\% & 1.11205e-15 & 6.81956e-16 & 4.95078e-16 \\
\quad 40000 & 11 & 1963 & \textbf{1.00926} ± 1.607\% & 3.94656e-17 & 3.33882e-17 & 2.47074e-17 \\
\quad 90000 & 11 & 2432 & \textbf{2.37411} ± 0.4164\% & 1.31937e-17 & 1.36944e-17 & 1.01987e-17 \\
PMHSS-GMRES & & & & & & & \\
\hline
\quad 10000 & 9 & 2204 & 0.345996 ± 0.02171\% & 1.96557e-14 & 1.37958e-14 & 1.00285e-14 \\
\quad 40000 & 10 & 3203 & 1.36888 ± 0.009568\% & 5.67274e-16 & 5.10832e-16 & 3.77199e-16 \\
\quad 90000 & 10 & 3995 & 3.27446 ± 0.01784\% & 1.92527e-16 & 1.97563e-16 & 1.46642e-16 \\
PRESB-GMRES & & & & & & & \\
\hline
\quad 10000 & 8 & 3464 & 0.258084 ± 0.008148\% & 3.35243e-15 & 4.61427e-15 & 2.46771e-15 \\
\quad 40000 & 8 & 5637 & 1.28754 ± 0.02749\% & 1.1019e-15 & 2.12654e-15 & 1.18669e-15 \\
\quad 90000 & 8 & 7030 & 3.05243 ± 0.02684\% & 3.76795e-16 & 8.44006e-16 & 4.7808e-16 \\
\hline \\
& & & \underbar{\textbf{Shifted-$\omega$ system}} \\
PMHSS\\
\hline 
\quad 10000 & 49 & 6281 & 1.085 ± 0.2419\% & 2.58145e-10 & 7.72275e-11 & 7.63703e-11 \\ 
\quad 40000 & 50 & 7126 & 3.157 ± 0.2576\% & 9.04324e-11 & 2.94851e-11 & 2.90853e-11 \\ 
\quad 90000 & 50 & 7187 & 6.058 ± 0.6971\% & 6.02718e-11 & 1.97633e-11 & 1.9518e-11 \\ 
GMRES & & & & & & &  \\
\hline
\quad 10000 & - & 212 & 1.677 ± 3.722\% & 7.91011e-11 & 2.32285e-09 & - \\ 
\quad 40000 & - & 232 & 5.214 ± 0.508\% & 3.87483e-11 & 1.35325e-09 & - \\ 
\quad 90000 & - & 233 & 10.05 ± 1.007\% & 2.6577e-11 & 9.42004e-10 & - \\ 
AA-PMHSS & & & & & & &  \\
\hline 
\quad 10000 & 18 & 2751 & 0.5119 ± 5.645\% & 1.61867e-11 & 1.2482e-11 & 1.03118e-11 \\ 
\quad 40000 & 21 & 3498 & 1.783 ± 1.185\% & 2.50785e-11 & 2.48957e-11 & 2.05558e-11 \\ 
\quad 90000 & 22 & 3648 & 3.639 ± 0.7017\% & 1.28299e-11 & 1.19634e-11 & 1.00378e-11 \\ 
PMHSS-GMRES & & & & & & &  \\
\hline
\quad 10000 & 18 & 4340 & 0.7649 ± 0.004377\% & 1.06948e-10 & 1.1078e-10 & 8.60833e-11 \\ 
\quad 40000 & 22 & 5572 & 2.461 ± 0.03121\% & 4.01809e-11 & 4.64125e-11 & 3.75068e-11 \\ 
\quad 90000 & 22 & 5906 & 5.013 ± 0.03158\% & 4.72434e-11 & 5.3194e-11 & 4.32722e-11 \\ 
PRESB-GMRES & & & & & & &  \\
\hline
\quad 10000 & 12 & 5616 & \textbf{0.4379} ± 0.009439\% & 5.23277e-11 & 7.00731e-11 & 4.42568e-11 \\ 
\quad 40000 & 12 & 6793 & \textbf{1.561 }± 0.02302\% & 5.95516e-11 & 7.75276e-11 & 4.89482e-11 \\ 
\quad 90000 & 12 & 6544 & \textbf{2.866 }± 0.03836\% & 6.26606e-12 & 7.18707e-11 & 4.42628e-11 \\ 
\hline \\
& & & \underbar{\textbf{Eq. of Motion}} \\
PMHSS \\
\hline 
\quad 10000 & 49 & 10358 & 1.786 ± 1.65\% & 1.25823e-10 & 7.30042e-11 & 7.09021e-11 \\ 
\quad 40000 & 51 & 20638 & 8.945 ± 1.766\% & 3.13699e-11 & 3.25443e-11 & 3.16235e-11 \\ 
\quad 90000 & 52 & 30185 & 24.9 ± 0.1285\% & 1.45818e-11 & 2.19632e-11 & 2.13458e-11 \\ 
GMRES & & & & & & &  \\
\hline
\quad 10000 & - & 280 & 3.057 ± 0.8814\% & 7.83319e-11 & 2.34518e-09 & - \\ 
\quad 40000 & - & 548 & 28.78 ± 1.626\% & 3.97398e-11 & 4.39533e-09 & - \\ 
\quad 90000 & - & 807 & 117.2 ± 1.591\% & 2.63963e-11 & 6.403e-09 & - \\ 
AA-PMHSS & & & & & & & \\
\hline 
\quad 10000 & 12 & 2732 & \textbf{0.5376} ± 3.711\% & 4.22335e-13 & 6.91976e-13 & 5.76201e-13 \\ 
\quad 40000 & 12 & 5369 & \textbf{2.679} ± 1.429\% & 5.5309e-12 & 4.84056e-12 & 4.36451e-12 \\ 
\quad 90000 & 12 & 7855 & \textbf{7.421} ± 0.6799\% & 5.67128e-12 & 4.05111e-12 & 3.65619e-12 \\ 
PMHSS-GMRES & & & & & & & \\
\hline
\quad 10000 & 11 & 4096 & 0.7076 ± 0.002372\% & 2.31419e-10 & 2.14003e-10 & 1.91861e-10 \\ 
\quad 40000 & 11 & 8092 & 3.542 ± 0.007765\% & 1.2527e-10 & 1.82148e-10 & 1.63178e-10 \\ 
\quad 90000 & 11 & 12055 & 9.983 ± 0.05302\% & 7.93328e-11 & 1.53629e-10 & 1.38611e-10 \\ 
PRESB-GMRES & & & & & & & \\
\hline
\quad 10000 & 11 & 8210 & 0.6628 ± 0.05686\% & 1.59327e-10 & 2.3776e-10 & 1.60705e-10 \\ 
\quad 40000 & 11 & 16261 & 3.716 ± 0.04519\% & 1.33754e-10 & 2.88983e-10 & 2.06747e-10 \\ 
\quad 90000 & 11 & 24320 & 10.34 ± 0.0956\% & 4.37123e-11 & 1.14813e-10 & 8.5615e-11 \\ 
\hline 
\end{tabular}
\end{scriptsize}
\end{table*}

\begin{table*}[]
\centering
\caption{The shifted-$\omega$ benchmark with CG for solving the inner system. The inner solvers' tolerance is set to $10^{-12}$, and the maximum number of inner iterations is limited to 50. The reference solution $x^*$ is calculated with Matlab's \texttt{backslash}.\\
} 
\label{table:perf_m50}
\begin{scriptsize}
\begin{tabular}{lccccccc}
Method/size & outer iter. & inner iter. & T ± $\sigma$ [s] & $||r||$ & $||x^*-x||$ & $||(A+B)^{-1}r||$ \\ 
\hline \\
& & & \underbar{\textbf{Shifted-$\omega$ system}} \\
PMHSS\\
\hline 
\quad 10000 & 49 & 2402 & 0.4153 ± 2.971\% & 2.59308e-10 & 7.84809e-11 & 7.76181e-11 \\ 
\quad 40000 & 49 & 2439 & 1.143 ± 0.6965\% & 1.2883e-10 & 3.98697e-11 & 3.93964e-11 \\ 
\quad 90000 & 49 & 2450 & 2.159 ± 1.086\% & 8.5479e-11 & 2.69666e-11 & 2.6633e-11 \\ 
AA-PMHSS & & & & & & &  \\
\hline 
\quad 10000 & 21 & 1049 & \textbf{0.2244} ± 3.837\% & 4.67345e-12 & 4.88027e-11 & 4.2427e-11 \\ 
\quad 40000 & 25 & 1248 & \textbf{0.8066} ± 1.546\% & 1.16763e-12 & 1.0343e-11 & 8.89348e-12 \\ 
\quad 90000 & 26 & 1299 & \textbf{1.78} ± 1.487\% & 1.12012e-12 & 7.91642e-12 & 6.65226e-12 \\ 
PMHSS-GMRES & & & & & & &  \\
\hline
\quad 10000 & 22 & 1394 & \text{stagnation}  & 4.37011e-06 & 7.87097e-05 & 6.59934e-05 \\ 
\quad 40000 & 26 & 1594 & \text{stagnation} & 2.76126e-06 & 3.71678e-05 & 3.17662e-05 \\ 
\quad 90000 & 26 & 1594 & \text{stagnation} & 2.71018e-06 & 2.47443e-05 & 2.07453e-05 \\ 
PRESB-GMRES & & & & & & &  \\
\hline
\quad 10000 & 16 & 2188 & \text{stagnation}  & 6.3969e-06 & 5.85338e-05 & 7.22871e-05 \\ 
\quad 40000 & 17 & 2293 & \text{stagnation}  & 4.05066e-06 & 3.94683e-05 & 6.08171e-05 \\ 
\quad 90000 & 17 & 2291 & \text{stagnation} & 2.91783e-06 & 3.06526e-05 & 4.62803e-05 \\ 
\hline \\
\end{tabular}
\end{scriptsize}
\end{table*}

\section{Related Work}
\label{rel_work}
Much effort has been allocated to solving complex linear systems and improving the conditioning of such systems as they naturally appear in physics related to waves, such as in electromagnetism~\cite{bai2012block}.  
A number of surveys on solvers and preconditioners for complex linear systems  can be found~\cite{axelsson2000real,axelsson2014comparison,wang2007short}. 

Several papers have been published on the HSS \cite{bai2007accelerated,bai2007convergence} family of solvers, including work on MHSS~\cite{bai2010modified} and the PMHSS with different choices of $V$ and $\alpha$ both as an independent solver and as a preconditioner for GMRES \cite{bai2011preconditioned}. The inner solver has been solved both with direct solvers and with iterative solvers, whereas the latter has been denoted with the prefix I, like IHSS. Some tested inner solvers are CG, Chebyshev iteration, and Algebraic Multigrid including some discussion on the tolerance of the inner solver \cite{axelsson2014comparison}. Recent work has suggested variants of the PMHSS \cite{zheng2022variant} and the use of PMHSS iteration preconditioners for Stokes control PDE constrained optimization problems \cite{cao2021pmhss}. 

Multiple works have noted or investigated the relationship between AA accelerated linear fixed-point iterations and GMRES, such as in \cite{walker2011anderson,doi:10.1137/S106482759426955X,fang2009two}. The acceleration of Jacobi iteration, another splitting method with Anderson acceleration, was developed in the paper~\cite{pratapa2016anderson}, which does not benefit from the accelerated solution of the preconditioner as Jacobi does not require a linear system to be solved.

\section{Discussion and Conclusions}
\label{conclusion}
This paper introduced a new technique for solving linear systems characterized by complex-valued matrices. Essentially, the new method is based on the PMHSS iteration and uses AA to accelerate the fixed-point iteration. 

We tested the new method against a number of benchmark problems: a time-harmonic eddy current simulation discretized with FEM and three benchmark problems found in the literature based in FD discretization. We can conclude that the combination of AA and PMHSS creates a fast and, for the tested problems, robust solver. Our method retains similar spectral properties and convergence as the PMHSS preconditioned GMRES method that has been proposed in previous work. However, it is significantly cheaper for systems that allow an iterative inner solver, which is necessary for large-scale simulations. The AA-PMHSS method does sacrifice some orthogonality in the Krylov space and speed of solving the minimization problem. However, these shortcomings are compensated by the faster evaluation of the precondition. It would be possible to achieve similar behavior with PMHSS-GMRES by allowing an adaptive tolerance for the inner solver. However, as our results show, the preconditioner becomes more expensive with the number of outer iterations as higher accuracy is required (at least the same accuracy as the error of the current solution). This is in contrast to AA-PMHSS, where the preconditioner becomes cheaper to evaluate throughout the outer iterations, due to being able to utilize the solution in the previous outer iteration as a starting point and not needing to fully resolve the inner iteration to enhance the solution. 

We also compared our method with the PRESB precondition, which has the advantage of forcing the spectrum to a smaller bound and guaranteeing real eigenvalues. In some of the tested cases, PRESB-GMRES does perform very well, owing to the reduced number of outer iterations required for convergence. For cases where the number of outer iterations is similar to the other methods, however, PRESB-GMRES falls behind since each outer iteration is more expensive to perform. PRESB requires two inner systems to be solved, both of which we cannot assume much about or speed up easily. It does also require the preconditioning system to be evaluated to a higher tolerance than the AA-PMHSS at each iteration to guarantee convergence to a sufficient stated tolerance. In the end, AA-PMHSS had the fastest wall-clock time in all but one test problem (shifted-$\omega$) when the inner solver was fully resolved and was the fastest in all test problems when capping the number of inner iterations to 50 (being 40\% faster than PRESB); only AA-PMHSS and PMHSS converged without stagnation. 

Many exciting research directions have unraveled from these examples, such as utilizing restarted methods (\textit{windowed} AA) to further negate the penalty from the more expensive and less robust minimization solution. There is also a need to formalize the relationship between AA and GMRES for the complex case. 

Furthermore, the idea of using Anderson acceleration to accelerate the convergence of fixed-point iterations in the context of splitting methods for linear systems is not limited to the PMHSS method alone. Rather, the approach could be used together with the more general HSS family of solvers or other methods based on matrix-splitting, which requires an iterative inner solver.

\section*{Acknowledgment}
The authors want to thank Martin Karp for reading and commenting on the manuscript.  
\bibliographystyle{siam}
\bibliography{main_siam.bib}

\begin{thebibliography}{10}

\bibitem{adams2007algebraic}
{\sc M.~F. Adams}, {\em Algebraic multigrid methods for direct frequency
  response analyses in solid mechanics}, Computational Mechanics, 39 (2007),
  pp.~497--507.

\bibitem{anderson1965iterative}
{\sc D.~G. Anderson}, {\em Iterative procedures for nonlinear integral
  equations}, Journal of the ACM (JACM), 12 (1965), pp.~547--560.

\bibitem{AXELSSON2020286}
{\sc O.~Axelsson}, {\em Optimality properties of a square block matrix
  preconditioner with applications}, Computers \& Mathematics with
  Applications, 80 (2020), pp.~286--294.
\newblock Numerical Methods for Scientific Computations and Advanced
  Applications II.

\bibitem{axelsson2000real}
{\sc O.~Axelsson and A.~Kucherov}, {\em Real valued iterative methods for
  solving complex symmetric linear systems}, Numerical linear algebra with
  applications, 7 (2000), pp.~197--218.

\bibitem{axelsson2014comparison}
{\sc O.~Axelsson, M.~Neytcheva, and B.~Ahmad}, {\em A comparison of iterative
  methods to solve complex valued linear algebraic systems}, Numerical
  Algorithms, 66 (2014), pp.~811--841.

\bibitem{bai2012block}
{\sc Z.-Z. Bai}, {\em Block alternating splitting implicit iteration methods
  for saddle-point problems from time-harmonic eddy current models}, Numerical
  Linear Algebra with Applications, 19 (2012), pp.~914--936.

\bibitem{bai2010modified}
{\sc Z.-Z. Bai, M.~Benzi, and F.~Chen}, {\em Modified {HSS} iteration methods
  for a class of complex symmetric linear systems}, Computing, 87 (2010),
  pp.~93--111.

\bibitem{bai2011preconditioned}
\leavevmode\vrule height 2pt depth -1.6pt width 23pt, {\em On preconditioned
  {MHSS} iteration methods for complex symmetric linear systems}, Numerical
  Algorithms, 56 (2011), pp.~297--317.

\bibitem{bai2007convergence}
{\sc Z.-Z. Bai, G.~Golub, and C.-K. Li}, {\em Convergence properties of
  preconditioned {H}ermitian and skew-{H}ermitian splitting methods for
  non-{H}ermitian positive semidefinite matrices}, Mathematics of Computation,
  76 (2007), pp.~287--298.

\bibitem{bai2007accelerated}
{\sc Z.-Z. Bai and G.~H. Golub}, {\em Accelerated {H}ermitian and
  skew-{H}ermitian splitting iteration methods for saddle-point problems}, IMA
  Journal of Numerical Analysis, 27 (2007), pp.~1--23.

\bibitem{baumann2018convergence}
{\sc M.~Baumann and M.~B. van Gijzen}, {\em Convergence and complexity study of
  {GMRES} variants for solving multi-frequency elastic wave propagation
  problems}, Journal of computational science, 26 (2018), pp.~285--293.

\bibitem{benzi2008block}
{\sc M.~Benzi and D.~Bertaccini}, {\em Block preconditioning of real-valued
  iterative algorithms for complex linear systems}, IMA Journal of Numerical
  Analysis, 28 (2008), pp.~598--618.

\bibitem{cao2021pmhss}
{\sc S.-M. Cao and Z.-Q. Wang}, {\em Pmhss iteration method and preconditioners
  for stokes control {PDE}-constrained optimization problems}, Numerical
  Algorithms, 87 (2021), pp.~365--380.

\bibitem{doi:10.1137/S106482759426955X}
{\sc N.~N. Carlson and K.~Miller}, {\em Design and {A}pplication of a
  {G}radient-weighted {M}oving {F}inite {E}lement {C}ode {I}: in {O}ne
  {D}imension}, SIAM Journal on Scientific Computing, 19 (1998), pp.~728--765.

\bibitem{de2021anderson}
{\sc H.~De~Sterck and Y.~He}, {\em Anderson acceleration as a {K}rylov method
  with application to asymptotic convergence analysis}, arXiv preprint
  arXiv:2109.14181,  (2021).

\bibitem{fang2009two}
{\sc H.-r. Fang and Y.~Saad}, {\em Two classes of multisecant methods for
  nonlinear acceleration}, Numerical linear algebra with applications, 16
  (2009), pp.~197--221.

\bibitem{feriani2000iterative}
{\sc A.~Feriani, F.~Perotti, and V.~Simoncini}, {\em Iterative system solvers
  for the frequency analysis of linear mechanical systems}, Computer Methods in
  Applied Mechanics and Engineering, 190 (2000), pp.~1719--1739.

\bibitem{li2018fast}
{\sc M.~Li, X.-M. Gu, C.~Huang, M.~Fei, and G.~Zhang}, {\em A fast linearized
  conservative finite element method for the strongly coupled nonlinear
  fractional {S}chr{\"o}dinger equations}, Journal of Computational Physics,
  358 (2018), pp.~256--282.

\bibitem{napov2012algebraic}
{\sc A.~Napov and Y.~Notay}, {\em An algebraic multigrid method with guaranteed
  convergence rate}, SIAM journal on scientific computing, 34 (2012),
  pp.~A1079--A1109.

\bibitem{ni2009anderson}
{\sc P.~Ni}, {\em Anderson acceleration of fixed-point iteration with
  applications to electronic structure computations}, PhD thesis, Worcester
  Polytechnic Institute Worcester, MA, USA, 2009.

\bibitem{notay2010aggregation}
{\sc Y.~Notay}, {\em An aggregation-based algebraic multigrid method},
  Electronic transactions on numerical analysis, 37 (2010), pp.~123--146.

\bibitem{notay2012aggregation}
\leavevmode\vrule height 2pt depth -1.6pt width 23pt, {\em Aggregation-based
  algebraic multigrid for convection-diffusion equations}, SIAM journal on
  scientific computing, 34 (2012), pp.~A2288--A2316.

\bibitem{pratapa2016anderson}
{\sc P.~P. Pratapa, P.~Suryanarayana, and J.~E. Pask}, {\em Anderson
  acceleration of the {J}acobi iterative method: {A}n efficient alternative to
  {K}rylov methods for large, sparse linear systems}, Journal of Computational
  Physics, 306 (2016), pp.~43--54.

\bibitem{rodriguez2010eddy}
{\sc A.~A. Rodr{\'\i}guez and A.~Valli}, {\em {E}ddy {C}urrent {A}pproximation
  of {M}axwell {E}quations: {T}heory, {A}lgorithms and {A}pplications}, vol.~4,
  Springer Science \& Business Media, 2010.

\bibitem{saad1986gmres}
{\sc Y.~Saad and M.~H. Schultz}, {\em {GMRES}: A generalized minimal residual
  algorithm for solving nonsymmetric linear systems}, SIAM Journal on
  scientific and statistical computing, 7 (1986), pp.~856--869.

\bibitem{walker1988implementation}
{\sc H.~F. Walker}, {\em Implementation of the {GMRES} method using
  {H}ouseholder transformations}, SIAM Journal on Scientific and Statistical
  Computing, 9 (1988), pp.~152--163.

\bibitem{walker2011anderson}
{\sc H.~F. Walker and P.~Ni}, {\em Anderson acceleration for fixed-point
  iterations}, SIAM Journal on Numerical Analysis, 49 (2011), pp.~1715--1735.

\bibitem{wang2007short}
{\sc Y.~Wang, J.~Lee, and J.~Zhang}, {\em A short survey on preconditioning
  techniques for large-scale dense complex linear systems in electromagnetics},
  International Journal of Computer Mathematics, 84 (2007), pp.~1211--1223.

\bibitem{wu2017modified}
{\sc S.-L. Wu and C.-X. Li}, {\em Modified complex-symmetric and
  skew-{H}ermitian splitting iteration method for a class of complex-symmetric
  indefinite linear systems}, Numerical Algorithms, 76 (2017), pp.~93--107.

\bibitem{zheng2022variant}
{\sc Z.~Zheng, M.-L. Zeng, and G.-F. Zhang}, {\em A variant of {PMHSS}
  iteration method for a class of complex symmetric indefinite linear systems},
  Numerical Algorithms, 91 (2022), pp.~283--300.

\end{thebibliography}
\end{document}